\documentclass[lettersize,journal]{IEEEtran}

\usepackage{cite}

\ifCLASSINFOpdf
   \usepackage[pdftex]{graphicx}
\else
   \usepackage[dvips]{graphicx}
\fi

\usepackage[cmex10]{amsmath}
\usepackage{amsthm}
\usepackage{amsfonts}
\usepackage{bm}

\usepackage{xcolor}
\usepackage{tikz}
\usetikzlibrary{calc}
\usetikzlibrary{shapes}
\usetikzlibrary{arrows}
\usetikzlibrary{positioning}
\usepackage{circuitikz}
\usepackage{enumerate}
\usepackage{enumitem}
\usepackage{algpseudocode}
\usepackage{algorithm}
\usepackage{empheq}
\usepackage{mathtools}
\usepackage{amssymb}

\allowdisplaybreaks

\theoremstyle{definition} 

\theoremstyle{plain} 
\newtheorem{proposition}{Proposition}
\newtheorem{theorem}{Theorem}

\theoremstyle{remark} 

\newcommand{\norm}[1]{\left\lVert#1\right\rVert}

\makeatletter
\let\old@ps@headings\ps@headings
\let\old@ps@IEEEtitlepagestyle\ps@IEEEtitlepagestyle
\def\psccfooter#1{%
    \def\ps@headings{%
        \old@ps@headings%
        \def\@oddfoot{\strut\hfill#1\hfill\strut}%
        \def\@evenfoot{\strut\hfill#1\hfill\strut}%
    }%
    \def\ps@IEEEtitlepagestyle{%
        \old@ps@IEEEtitlepagestyle%
        \def\@oddfoot{\strut\hfill#1\hfill\strut}%
        \def\@evenfoot{\strut\hfill#1\hfill\strut}%
    }%
    \ps@headings%
}
\makeatother

\begin{document}

\title{Second-Order Cone Valid Inequalities for Set-Point Tracking with Energy Storage}

\author{Juan M. Morales
\thanks{\noindent
J. M. Morales is with the research group OASYS, University of Málaga, Málaga 29071, Spain. E-mail: juan.morales@uma.es.\\
This work was supported by the Consejería de Universidad, Investigación e Innovación de la Junta de Andalucía through project DGP-PIDI-2024-00851. The author thankfully acknowledges the computer resources, technical expertise, and assistance provided by the SCBI (Supercomputing and Bioinformatics) center of the University of M\'alaga.}
}

\maketitle

\begin{abstract}
Convex relaxations of energy storage models drop the complementarity constraints that forbid simultaneous charging and discharging. In the set-point tracking problem, the relaxation can then charge and discharge simultaneously at no cost, as the tracking error depends only on the net power of the storage. We show that the convex hull of the single-period epigraph of the tracking error, under complementarity, is the epigraph of a parabolic cylinder. The resulting second-order cone inequality adds to the tracking error a penalty on the product of the charging and discharging powers and is the tightest convex member of the pencil of quadrics through the two branches of the epigraph. With power limits and the mode variable, the hull is given by the perspective reformulation of the two modes, which dominates the cylinder. In 20\,000 instances, the cylinder reduces the hours with simultaneous charging and discharging by more than 99\% and yields an optimal schedule in 98.8\% of them.
\end{abstract}


\section{Introduction}\label{sec:intro}
Storage models need complementarity constraints, written with bilinear equalities or binary variables, to prevent charging and discharging in the same period; their convex relaxations drop them and may return infeasible schedules~\cite{arroyo2020use}. Valid inequalities that strengthen the relaxation reduce this effect: Pozo~\cite{pozo2023convex} describes the convex hull of the single-period storage set, and Elgersma \emph{et al.}~\cite{elgersma2024tight} extend it to reserve provision. These inequalities act on the feasible set and leave the objective function untouched. Al Taha and Bitar~\cite{taha2024lossy} instead eliminate the power variables through a piecewise-affine bijection with the state-of-charge profiles, whose feasible set is convex; the nonconvexity remains in the mapping, so the power variables may appear only in the objective function, which stays convex if it is nondecreasing in the charging power. In the set-point tracking problem, the relaxation is weak precisely because of its objective: the tracking error depends only on the net power, so charging and discharging simultaneously is free~\cite{pozo2023convex}; it also violates the condition of~\cite{taha2024lossy} whenever the set-point asks the storage to charge.

This letter shows that the convex hull of the single-period epigraph of the tracking error under complementarity is the epigraph of a parabolic cylinder, which yields a second-order cone (SOC) inequality that penalizes simultaneous charging and discharging (Section~\ref{sec:hull}). We identify this cylinder as the tightest convex member of the pencil of quadrics through the epigraph, extend the result to power limits through the perspective reformulation of the operating modes, and show that the cut nearly eliminates simultaneous charging and discharging in 20\,000 instances (Section~\ref{sec:case_study}).

\section{Set-Point Tracking Problem}\label{sec:Formulation}
Let $p_t^d$, $p_t^c$, and $s_t$ denote the discharging power, the charging power, and the state of charge of the storage in period $t \in \{1, \ldots, T\}$. Its operational set is
\begin{subequations}\label{eq:Pset}
  \begin{align}
      \mathcal{P} = \Big\{&(\boldsymbol{p}^d, \boldsymbol{p}^c, \boldsymbol{s}) \in \mathbb{R}^{3T}_{\geq 0}: \ \boldsymbol{p}^d \leq \boldsymbol{1}\, \overline{P}^d, \ \boldsymbol{p}^c \leq \boldsymbol{1}\, \overline{P}^c, \label{eq:STP_max_power}\\
&\quad \boldsymbol{1}\, \underline{S} \leq \boldsymbol{s} \leq \boldsymbol{1}\, \overline{S}, \label{eq:STP_max_level}\\
& \quad s_{t} = s_{t-1} + \Delta\big(\eta_{c} p_t^c - p_t^d/\eta_d\big), \ \ \forall t, \label{eq:STP_dynamic_state}\\
& \quad  p_t^d\, p_t^{c} = 0, \ \ \forall t\Big\},  \label{eq:STP_comp}
  \end{align}
\end{subequations}
where $\overline{P}^d$ and $\overline{P}^c$ are the power limits, $[\underline{S}, \overline{S}]$ is the state-of-charge range, $\eta_d, \eta_c \in (0,1]$ are the efficiencies, $\Delta$ is the period length, and $s_0$ is given. The complementarity constraints~\eqref{eq:STP_comp} can be written with a binary variable $u_t$ per period and the power limits $p_t^d \leq \overline{P}^d (1-u_t)$ and $p_t^c \leq \overline{P}^c u_t$~\cite{pozo2023convex}.

The set-point tracking problem (STP) seeks the operation of the storage that best follows a power signal $p_t^s$~\cite{pozo2023convex}:
\begin{equation}\label{eq:STP}
 \min_{(\boldsymbol{p}^d, \boldsymbol{p}^c, \boldsymbol{s}) \in \mathcal{P}}  \quad \sum_{t=1}^T \big(p_t^d-p_t^c-p_t^s\big)^2.
\end{equation}
Its natural relaxation drops~\eqref{eq:STP_comp} (or relaxes $u_t$ to $[0,1]$). As the tracking error depends only on the net power, the relaxation can increase both $p_t^d$ and $p_t^c$ at no cost, and often does so to meet the state-of-charge limits~\cite{pozo2023convex}. To strengthen the objective, we write~\eqref{eq:STP} in epigraph form, with one variable $z_t$ per period:
\begin{equation}\label{eq:STP_epi}
\min_{(\boldsymbol{p}^d, \boldsymbol{p}^c, \boldsymbol{s}) \in \mathcal{P},\, \boldsymbol{z}} \Big\{ \sum_{t=1}^T z_t : (z_t, p_t^d, p_t^c) \in \mathcal{E}_t \ \ \forall t \Big\},
\end{equation}
where
\begin{equation}\label{eq:E}
\mathcal{E}_t := \big\{(z, p^d, p^c) \in \mathbb{R} \times \mathbb{R}^2_{\geq 0} : z \geq (p^d - p^c - p_t^s)^2, \ p^d p^c = 0\big\}.
\end{equation}
As~\eqref{eq:STP_epi} is linear in $z_t$, the tightest convex relaxation that can be built from each set $\mathcal{E}_t$ alone is its convex hull.

\section{The Convex Hull of the Epigraph}\label{sec:hull}
\begin{figure}
\centering
\begin{tikzpicture}[x={(0.80cm,-0.16cm)}, y={(-0.58cm,-0.26cm)}, z={(0cm,0.25cm)}, line cap=round, line join=round]
  \foreach \r in {0.3,0.6,0.9,1.2,1.5} {
    \draw[black!55, densely dashed, thin, domain=0:1, samples=30, smooth] plot ({\x*\r}, {(1-\x)*\r}, {(\r*(2*\x-1)-1)^2});
  }
  \foreach \ra/\rb in {0/0.3,0.3/0.6,0.6/0.9,0.9/1.2,1.2/1.5} {
    \fill[violet, opacity=0.18] (\ra,0,{(\ra-1)^2}) -- (\rb,0,{(\rb-1)^2}) -- (0,\rb,{(\rb+1)^2}) -- (0,\ra,{(\ra+1)^2}) -- cycle;
    \draw[violet, thin] (\rb,0,{(\rb-1)^2}) -- (0,\rb,{(\rb+1)^2});
  }
  \draw[very thick, domain=0:1.5, samples=40, smooth] plot (\x, 0, {(\x-1)^2});
  \draw[very thick, domain=0:1.5, samples=40, smooth] plot (0, \x, {(\x+1)^2});
  \draw[->] (0,0,0) -- (1.8,0,0) node[right, font=\scriptsize]{$p^d$};
  \draw[->] (0,0,0) -- (0,1.8,0) node[left, font=\scriptsize]{$p^c$};
  \draw[->] (0,0,0) -- (0,0,7) node[above, font=\scriptsize]{$z$};
  \node[font=\scriptsize] at (0.9,0.9,-2.2) {(a)};
\end{tikzpicture}\hfill
\begin{tikzpicture}[x=1.6cm, y=0.25cm]
  \fill[violet!12, domain=0:1, samples=60] plot (\x, {6.25 - 6*\x}) -- plot[domain=1:0, samples=60] (\x, {(1.5*(2*\x-1)-1)^2}) -- cycle;
  \draw[->] (0,0) -- (1.12,0) node[right, font=\scriptsize]{$\theta$};
  \draw[->] (0,0) -- (0,6.9) node[above, font=\scriptsize]{$z$};
  \foreach \v in {2,4,6} \draw (0,\v) -- (-0.03,\v) node[left, font=\scriptsize]{\v};
  \draw (0,0) -- (0,-0.15) node[below, font=\scriptsize]{0};
  \draw (1,0) -- (1,-0.15) node[below, font=\scriptsize]{1};
  \draw[violet, thick] (0,6.25) -- (1,0.25);
  \draw[black, thin, domain=0:1, samples=60, smooth] plot (\x, {2.25*(\x*\x+(1-\x)*(1-\x)) - 3*(2*\x-1) + 1});
  \draw[black!55, densely dashed, thick, domain=0:1, samples=60, smooth] plot (\x, {(1.5*(2*\x-1)-1)^2});
  \fill (0,6.25) circle (1.3pt); \fill (1,0.25) circle (1.3pt);
  \draw[violet, thin] (0.5,3.25) -- (0.62,5.6) node[right, font=\scriptsize, inner sep=1pt]{$\lambda = 2$};
  \draw[thin] (0.55,1.84) -- (0.74,4.2) node[right, font=\scriptsize, inner sep=1pt]{$0$};
  \draw[black!55, thin] (0.62,0.41) -- (0.86,2.8) node[right, font=\scriptsize, black, inner sep=1pt]{$-2$};
  \node[font=\scriptsize] at (0.5,-2.4) {(b)};
\end{tikzpicture}\hfill
\begin{tikzpicture}[x=1.15cm, y=0.98cm, line cap=round, line join=round]
  \fill[violet!12] (0,0) -- (1,0) -- (0,1) -- cycle;
  \fill[black!10] (1,0) -- (0,1.6) -- (0,1) -- cycle;
  \foreach \r in {0.25,0.5,0.75} \draw[violet, thin] (\r,0) -- (0,\r);
  \draw[violet, thick] (1,0) -- (0,1);
  \foreach \b in {1.2,1.4} \draw[thin] (1,0) -- (0,\b);
  \draw[thick] (1,0) -- (0,1.6);
  \fill (1,0) circle (1.3pt);
  \draw[->] (0,0) -- (1.3,0) node[right, font=\scriptsize, inner sep=1pt]{$p^d$};
  \draw[->] (0,0) -- (0,1.85) node[above, font=\scriptsize]{$p^c$};
  \node[below, font=\scriptsize] at (1,0) {$\overline{P}^d$};
  \node[left, font=\scriptsize, inner sep=1.5pt] at (0,1) {$\overline{P}^d$};
  \node[left, font=\scriptsize, inner sep=1.5pt] at (0,1.6) {$\overline{P}^c$};
  \node[font=\scriptsize] at (0.6,-0.62) {(c)};
\end{tikzpicture}
\caption{Tracking error for $p^s = 1$. (a) Branches of $\mathcal{E}$ (thick), rulings of the cylinder~\eqref{eq:cylinder} (solid), and natural relaxation along the segments $p^d + p^c = r$ (dashed). (b) Members $\lambda = 2$ (cylinder), $0$ (paraboloid), and $-2$ (natural relaxation) of the pencil~\eqref{eq:pencil} along $p^d + p^c = 1.5$, with $\theta = p^d/(p^d+p^c)$; the SOC cut removes the shaded area; dots: points of $\mathcal{E}$. (c) Top view of the hull with power limits $\overline{P}^d < \overline{P}^c$: rulings of the cylinder and, where $p^d + p^c > \overline{P}^d$, of the cone~\eqref{eq:cone} with apex at the full-discharge point (dot).}
\label{fig:cylinder}
\end{figure}
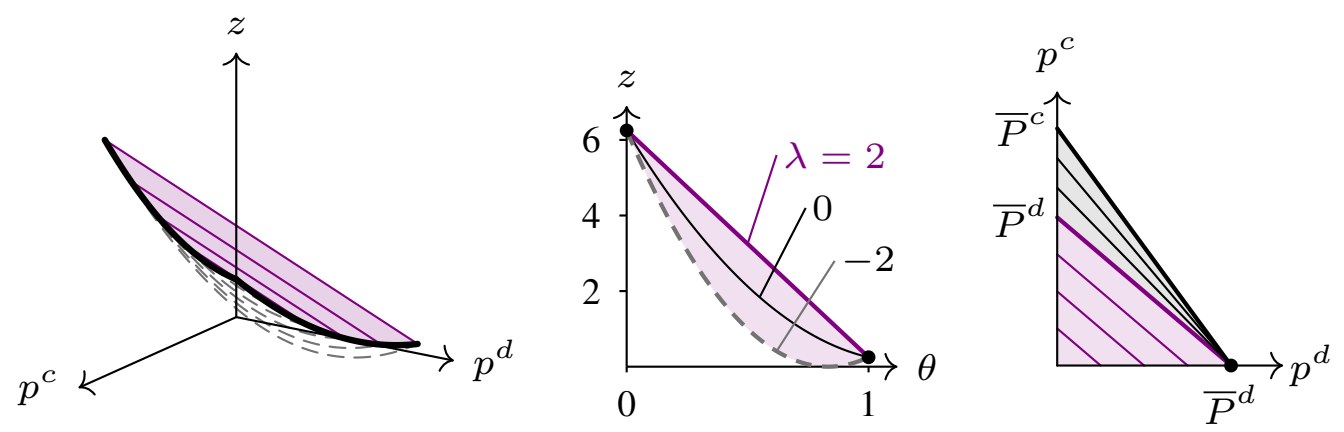

We omit the subscript $t$ and write $p^s$ for $p_t^s$. The set $\mathcal{E}$ is the union of two convex sets, the epigraphs of the pure-discharge parabola $z = (p^d - p^s)^2$ on the plane $p^c = 0$ and of the pure-charge parabola $z = (p^c + p^s)^2$ on the plane $p^d = 0$; these parabolas are the two branches of $\mathcal{E}$ (thick curves in Fig.~\ref{fig:cylinder}a). Its convex hull is the epigraph of the parabolic cylinder
\begin{equation}\label{eq:cylinder}
g_{\mathcal{C}}(p^d, p^c) := (p^d + p^c)^2 - 2p^s(p^d - p^c) + (p^s)^2.
\end{equation}

\begin{theorem}\label{Thm:convex_hull}
$\operatorname{conv}(\mathcal{E}) = \mathcal{K} := \big\{(z, p^d, p^c) \in \mathbb{R} \times \mathbb{R}^2_{\geq 0} : z \geq g_{\mathcal{C}}(p^d, p^c)\big\}$.
\end{theorem}

\begin{proof}
The set $\mathcal{K}$ is convex, because $g_{\mathcal{C}}$ is a convex quadratic function. On $\mathcal{E}$, $p^d p^c = 0$, so $(p^d + p^c)^2 = (p^d - p^c)^2$ and $g_{\mathcal{C}}(p^d,p^c) = (p^d - p^c - p^s)^2 \leq z$. Hence $\mathcal{E} \subseteq \mathcal{K}$ and $\operatorname{conv}(\mathcal{E}) \subseteq \mathcal{K}$. Conversely, let $(z, p^d, p^c) \in \mathcal{K}$, with throughput $r := p^d + p^c$ and slack $\delta := z - g_{\mathcal{C}}(p^d, p^c) \geq 0$. If $r = 0$, the point belongs to $\mathcal{E}$. Otherwise, let $\theta := p^d / r \in [0,1]$ and take the points that discharge or charge at power $r$, with their tracking errors raised by $\delta$,
\[
A := \big((r - p^s)^2 + \delta, r, 0\big), \quad B := \big((r + p^s)^2 + \delta, 0, r\big).
\]
Both belong to $\mathcal{E}$, as $\delta \geq 0$, and $\theta A + (1-\theta) B = (\zeta + \delta, p^d, p^c)$ with
\begin{align*}
\zeta &= \theta (r - p^s)^2 + (1-\theta)(r + p^s)^2 \\
&= r^2 - 2 p^s r (2\theta - 1) + (p^s)^2 = g_{\mathcal{C}}(p^d, p^c),
\end{align*}
because $r = p^d + p^c$ and $r(2\theta - 1) = p^d - p^c$. As $\zeta + \delta = z$, this convex combination of two points of $\mathcal{E}$ is $(z, p^d, p^c)$.
\end{proof}

Thus, along each segment $p^d + p^c = r$, the hull is the chord between the pure-discharge and pure-charge points with throughput $r$ (Fig.~\ref{fig:cylinder}a). The SOC cut $z_t \geq g_{\mathcal{C}}(p_t^d, p_t^c)$ can be written as
\begin{equation}\label{eq:SOC_cylinder}
z_t \geq \big(p_t^d - p_t^c - p_t^s\big)^2 + 4\, p_t^d\, p_t^c,
\end{equation}
that is, the tracking error plus a penalty on simultaneous charging and discharging, which vanishes on $\mathcal{E}_t$. It is a convex quadratic, hence SOC-representable, constraint that does not involve the binary variables, so it can be added to any relaxation of~$\mathcal{P}$. No epigraph variables are needed: replacing the objective of~\eqref{eq:STP} with $\sum_t \big[(p_t^d - p_t^c - p_t^s)^2 + 4 p_t^d p_t^c\big]$ yields the same relaxation as a convex quadratic program.

\paragraph*{Is the cylinder unique?} The cylinder belongs to the pencil of quadrics that contain the two branches of $\mathcal{E}$, which, following~\cite[Thm.~4.1]{belotti2013families}, we write as $z = g_\lambda(p^d, p^c)$ with
\begin{equation}\label{eq:pencil}
g_\lambda(p^d, p^c) := (p^d - p^c - p^s)^2 + (2 + \lambda)\, p^d p^c, \qquad \lambda \in \mathbb{R}.
\end{equation}
All members coincide on the planes $p^d = 0$ and $p^c = 0$, and the Hessian of $g_\lambda$, with diagonal entries $2$ and off-diagonal entries $\lambda$, is positive semidefinite if and only if $|\lambda| \leq 2$ and singular if and only if $\lambda = \pm 2$, the two cylinders of the pencil. For $\lambda = 2$, $g_2 = g_{\mathcal{C}}$ is the cylinder of Theorem~\ref{Thm:convex_hull}, with rulings parallel to $(-4p^s, 1, -1)$ in the space $(z, p^d, p^c)$; for $\lambda = -2$, $g_{-2}$ is the tracking error itself, i.e., the objective of the natural relaxation, with rulings parallel to $(0, 1, 1)$. Since $p^d p^c \geq 0$, $g_\lambda$ is non-decreasing in $\lambda$, so the cylinder $\lambda = 2$ is the tightest convex member and the cylinder $\lambda = -2$ the weakest (Fig.~\ref{fig:cylinder}b). Hence the cylinder is not unique, but it is the only one that yields the convex hull.

\subsection{Power Limits: Perspective Reformulation}\label{sec:perspective}
Theorem~\ref{Thm:convex_hull} ignores the power limits. With them, the cylinder remains valid but need not give the hull. In formulations that keep the mode variable $u$, the single-period set is
\begin{multline}\label{eq:M}
\mathcal{M} := \big\{(z, p^d, p^c, u) \in \mathbb{R} \times \mathbb{R}^2_{\geq 0} \times \{0,1\} : p^d \leq \overline{P}^d (1-u), \\
p^c \leq \overline{P}^c u, \ z \geq (p^d - p^c - p^s)^2 \big\},
\end{multline}
whose projection onto $(z, p^d, p^c)$ is $\mathcal{E}$ with the power limits. In each mode, the tracking error is a convex function of one power, $(p^d - p^s)^2$ or $(p^c + p^s)^2$, and the hull of such a disjunction follows from their perspectives~\cite{ceria1999convex}. Let
\begin{equation}\label{eq:persp}
h(p^d, p^c, u) := \frac{(p^d)^2}{1-u} + \frac{(p^c)^2}{u} - 2p^s(p^d - p^c) + (p^s)^2,
\end{equation}
with the convention $0^2/0 := 0$.

\begin{proposition}\label{prop:perspective}
$\operatorname{conv}(\mathcal{M}) = \mathcal{F} := \big\{(z, p^d, p^c, u) : 0 \leq p^d \leq \overline{P}^d (1-u), \ 0 \leq p^c \leq \overline{P}^c u, \ z \geq h(p^d, p^c, u)\big\}$. With two auxiliary variables $w^d, w^c \geq 0$, the inequality $z \geq h$ is equivalent to the SOC constraints
\begin{subequations}\label{eq:persp_soc}
\begin{align}
& z \geq w^d + w^c - 2p^s(p^d - p^c) + (p^s)^2, \label{eq:persp_soc_z}\\
& \norm{\big(2p^d,\ w^d - 1 + u\big)} \leq w^d + 1 - u, \label{eq:persp_soc_d}\\
& \norm{\big(2p^c,\ w^c - u\big)} \leq w^c + u, \label{eq:persp_soc_c}
\end{align}
\end{subequations}
where~\eqref{eq:persp_soc_d} and~\eqref{eq:persp_soc_c} are the rotated cones $(p^d)^2 \leq w^d (1-u)$ and $(p^c)^2 \leq w^c u$, as $\norm{(2x, w - y)} \leq w + y \Leftrightarrow x^2 \leq wy,\ w, y \geq 0$.
\end{proposition}

\begin{proof}
The set $\mathcal{F}$ is convex, as the projection of a set defined by linear and SOC constraints, and contains $\mathcal{M}$: if $u = 0$, then $p^c = 0$ and $h = (p^d - p^s)^2 \leq z$; if $u = 1$, then $p^d = 0$ and $h = (p^c + p^s)^2 \leq z$. Hence $\operatorname{conv}(\mathcal{M}) \subseteq \mathcal{F}$. Conversely, let $(z, p^d, p^c, u) \in \mathcal{F}$, with slack $\delta := z - h(p^d, p^c, u) \geq 0$. Its bounds give $u \in [0,1]$; if $u \in \{0,1\}$, it is in $\mathcal{M}$ by the cases above. Otherwise, $a := p^d/(1-u) \in [0, \overline{P}^d]$ and $b := p^c/u \in [0, \overline{P}^c]$, so $A := ((a - p^s)^2 + \delta, a, 0, 0)$ and $B := ((b + p^s)^2 + \delta, 0, b, 1)$ belong to $\mathcal{M}$, and $(1-u)A + uB = (\zeta + \delta, p^d, p^c, u)$ with $\zeta = (1-u)(a-p^s)^2 + u(b+p^s)^2 = (1-u)a^2 + ub^2 - 2p^s\big((1-u)a - ub\big) + (p^s)^2 = h(p^d, p^c, u)$, because $(1-u)a = p^d$ and $ub = p^c$. As $\zeta + \delta = z$, this convex combination of two points of $\mathcal{M}$ is $(z, p^d, p^c, u)$.
\end{proof}

Physically, a fractional $u$ splits the period: the storage discharges at power $a = p^d/(1-u)$ during a fraction $1 - u$ and charges at power $b = p^c/u$ during the rest, and $h$ is the time average of the two tracking errors, whereas the natural relaxation takes the error of the average net power.

The perspective inequality dominates the cylinder: by the Cauchy--Schwarz inequality, $\big((p^d)^2/(1-u) + (p^c)^2/u\big)\big((1-u) + u\big) \geq (p^d + p^c)^2$, so $h \geq g_{\mathcal{C}}$, with equality iff $a = b$, i.e., $u = u^\star := p^c/(p^d + p^c)$. The cylinder is thus the perspective inequality for the best time split, where both modes run at power $r = p^d + p^c$, as $A$ and $B$ in the proof of Theorem~\ref{Thm:convex_hull}.

The limits $a \leq \overline{P}^d$ and $b \leq \overline{P}^c$ confine $u$ to the interval $[p^c/\overline{P}^c,\ 1 - p^d/\overline{P}^d]$, which requires $p^d/\overline{P}^d + p^c/\overline{P}^c \leq 1$. As projection and convex hull commute, minimizing $h$ over this interval gives the hull of $\mathcal{E}$ with power limits (Fig.~\ref{fig:cylinder}c). Where $r \leq \min\{\overline{P}^d, \overline{P}^c\}$, the split $u^\star$ is feasible and the hull is the cylinder. Where, say, $r > \overline{P}^d$, it would discharge above $\overline{P}^d$ and, $h$ being convex in $u$, the minimum is attained at the closest feasible split, $1 - u = p^d/\overline{P}^d$, which discharges at full power, $a = \overline{P}^d$. Substituting it in~\eqref{eq:persp},
\begin{equation}\label{eq:cone}
z \geq \overline{P}^d p^d + \frac{\overline{P}^d (p^c)^2}{\overline{P}^d - p^d} - 2p^s(p^d - p^c) + (p^s)^2.
\end{equation}
These points lie on segments that join the full-discharge point $((\overline{P}^d - p^s)^2, \overline{P}^d, 0)$ with the pure-charge points, so~\eqref{eq:cone} is a cone with apex at the former: with $\ell$ denoting its affine terms, it reads $(z - \ell)(\overline{P}^d - p^d) \geq \overline{P}^d (p^c)^2$.

If $u_t$ appears only in the power limits, as in the natural relaxation, the optimizer chooses it to minimize $h$, and~\eqref{eq:persp_soc} yields the hull of Fig.~\ref{fig:cylinder}c: it tightens the cylinder only where $p^d + p^c > \min\{\overline{P}^d, \overline{P}^c\}$, and nowhere if $\overline{P}^d = \overline{P}^c$, since then $p^d/\overline{P}^d + p^c/\overline{P}^c \leq 1$ implies $p^d + p^c \leq \overline{P}^d$. They pay off when $u_t$ is coupled with other constraints (e.g., mode-dependent costs or reserves), at the price of two rotated cones per period; the cylinder needs neither $u_t$ nor auxiliary variables.

\section{Case Study}\label{sec:case_study}
As in~\cite{pozo2023convex}, a household with a 35-kW photovoltaic (PV) system uses a battery to follow its net demand over $T = 24$ hours; the first 200 daily PV profiles and the 100 battery configurations of~\cite{PozoData} yield 20\,000 instances. We solve the MIQP that models~\eqref{eq:STP_comp} with binary variables; the quadratic relaxation based on the HCH model of~\cite{pozo2023convex} (HCH-QP); HCH-QP with the cylinder as convex QP objective (HCH-SOC); and HCH-QP with~\eqref{eq:persp_soc} (HCH-PER), with Gurobi 13.0.2 (MOSEK 10.2.0 for HCH-PER) through AMPL, on a 2.6~GHz CPU with 2~GB RAM, one thread and default settings; times are solver CPU times reported by AMPL. See~\cite{mygithub} for data and code.

Table~\ref{tab:results_STP} reports, as in~\cite{pozo2023convex}, the hours with $p_t^c p_t^d > 10^{-4}$~kW$^2$, the average $\boldsymbol{p}^c \cdot \boldsymbol{p}^d$ and gap to the MIQP optimum, and the average time and its reduction with respect to the MIQP.
\begin{table}
\renewcommand{\arraystretch}{1.1}
\caption{Results for the set-point tracking problem (20\,000 instances)}
\label{tab:results_STP}
\centering
\setlength{\tabcolsep}{3pt}
\begin{tabular}{lccccc}
\hline
 & $\#[p_{t}^c p_{t}^d > 10^{-4}]$ & $\boldsymbol{p}^c \cdot \boldsymbol{p}^d$ & Gap & Time & Time red. \\
 & (\%) & (kW$^2$) & (\%) & (ms) & (\%) \\
\hline
MIQP & 0 & 0 & -- & 33.63 & -- \\
HCH-QP & 15.87 & 104.44 & $-1.34$ & 14.86 & 55.82 \\
HCH-SOC & 0.08 & 0.04 & $-0.0004$ & 14.66 & 56.40 \\
HCH-PER & 0.08 & 0.04 & $-0.0004$ & 20.39 & 39.36 \\
\hline
\end{tabular}
\end{table}
The cylinder removes 99.5\% of the hours with simultaneous charging and discharging and 99.96\% of $\boldsymbol{p}^c \cdot \boldsymbol{p}^d$, and nearly closes the gap of HCH-QP. In 19\,750 instances (98.8\%), HCH-SOC returns a schedule without simultaneous charging and discharging. Such a schedule is feasible for~\eqref{eq:STP} and, as the penalty in~\eqref{eq:SOC_cylinder} vanishes, its tracking error equals the lower bound given by the relaxation, so it is optimal; the MIQP confirms it (relative tolerance $10^{-4}$). The other 250 involve 11 small batteries with low efficiencies, which suggests a coupling of periods through the state-of-charge limits that a single-period hull cannot capture. HCH-PER gives the same bounds at a higher cost: $\overline{P}^d \neq \overline{P}^c$ in 99 batteries, but no hour with simultaneous charging and discharging reaches $p_t^d + p_t^c > \min\{\overline{P}^d, \overline{P}^c\}$, where the perspective is tighter (Section~\ref{sec:perspective}). HCH-SOC is as fast as HCH-QP, whose optima are not unique: the tracking error fixes only the net power $p_t^d - p_t^c$, so every feasible split of it into charging and discharging is optimal, and the solver picks one.

\section{Conclusions}\label{sec:conclusions}
The convex hull of the single-period epigraph of the tracking error under complementarity is the epigraph of a parabolic cylinder, the tightest convex quadric through its two branches. The resulting SOC inequality penalizes simultaneous charging and discharging without binary variables and nearly eliminates it. With power limits and the mode variable, the perspective reformulation gives the single-period hull. Future work will address several units and multi-period hulls.

\bibliographystyle{IEEEtran}
\bibliography{references_letter}

\end{document}